\def \version {2013--09--18}

\documentclass[12pt]{article}
\usepackage{amssymb}
\def \qed {\hfill $\boxempty$}

\newtheorem{Theorem}{Theorem}
\newtheorem{lem}[Theorem]{Lemma}
\newtheorem{defi}[Theorem]{Definition}
\newtheorem{crl}[Theorem]{Corollary}
\newtheorem{prp}[Theorem]{Proposition}
\newtheorem{prm}[Theorem]{Problem}
\newtheorem{cnj}[Theorem]{Conjecture}
\newtheorem{rmk}[Theorem]{Remark}
\newtheorem{xmp}[Theorem]{Example}
\newtheorem{obs}[Theorem]{Observation}

\def \bo {\begin{obs} \ }
\def \eo {\end{obs}}

\def \bp {\begin{prp} \ }
\def \ep {\end{prp}}

\def \bpm {\begin{prm} \ }
\def \epm {\end{prm}}

\def \bc {\begin{crl} \ }
\def \ec {\end{crl}}

\def \bcj {\begin{cnj} \ }
\def \ecj {\end{cnj}}

\def \bpm {\begin{prm} \ }
\def \epm {\end{prm}}

\def \thm {\begin{Theorem} \ }
\def \ethm {\end{Theorem}}

\def \bl {\begin{lem} \ }
\def \el {\end{lem}}

\def \bd {\begin{defi} \ \rm }
\def \ed {\end{defi}}

\def \brm {\begin{rmk} \ }
\def \erm {\end{rmk}}

\def \bxm {\begin{xmp} \ \rm }
\def \exm {\end{xmp}}

\def \nmr {\begin{enumerate}}
\def \enmr {\end{enumerate}}

\def \tmz {\begin{itemize}}
\def \etmz {\end{itemize}}

\def \smin {\setminus}
\def \nin {\noindent}
\def \bsk {\bigskip}
\def \msk {\medskip}

\def \pf {\nin{\bf Proof.} \ }

\def \qed {\hfill $\Box$}

\def \cH {{\cal H}}
\def \exsk {\mbox{\rm ex}(n,\cH^{(r)}(k,q))}
\def \ex {\mbox{\rm ex} }

\def \hnp {H^{(r)}_{n,p}}
\def \bbE {{\mathbb{E}}}

\def\cH{{\cal H}}

\def\cF{{\cal F}}
\def\cG{{\cal G}}

\def\cO{{\cal O}}

\newtheorem{claim}{Claim}

\def \bcl {\begin{claim} \ }
\def \ecl {\end{claim}}

\newtheorem{con}{Condition}

\def \bcon {\begin{con} \ \rm }
\def \econ {\end{con}}

\def \dia {\hfill $\Diamond$}

\begin{document}

\title{\vspace{-11ex}~~\\ Tur\'an numbers and batch codes\thanks{
  ~Research supported in part by the Hungarian Scientific Research
  Fund, OTKA grant T-81493, moreover
 by the European Union and Hungary, co-financed
 by the European Social Fund through the project T\'AMOP-4.2.2.C-11/1/KONV-2012-0004 -- National Research Center
 for Development and Market Introduction of Advanced Information and Communication
 Technologies.}}
\author{Csilla Bujt\'as~$^1$
    \qquad
        Zsolt Tuza~$^{1,2}$\\
\normalsize $^1$~Department of Computer Science and Systems
 Technology \\
  \normalsize University of Pannonia \\
\normalsize  Veszpr\'em,  Hungary \\
\normalsize $^2$~Alfr\'ed R\'enyi Institute of Mathematics \\
       \normalsize Hungarian Academy of Sciences \\
\normalsize  Budapest, Hungary
 }
\date{\footnotesize Latest update  on \version}
\maketitle

\begin{abstract}

Combinatorial batch codes provide a tool for distributed
 data storage, with the feature of keeping privacy during
 information retrieval.
Recently, Balachandran and Bhattacharya observed that the problem
 of constructing such uniform codes in an economic way can be
 formulated as a Tur\'an-type question on hypergraphs.
Here we establish general lower and upper bounds for this
 extremal problem, and also for its generalization where
 the forbidden family consists of those $r$-uniform hypergraphs
 $H$ which satisfy the condition $k\ge |E(H)|> |V(H)|+q$
  (for $k>q+r$ and $q> -r$ fixed).
We also prove that, in the given range of parameters,
 the considered Tur\'an function is asymptotically
 equal to the one restricted to $|E(H)|=k$,
  studied by Brown, Erd\H os and T.~S\'os.
Both families contain some $r$-partite members
 --- often called the `degenerate case',
 characterized by the equality $\lim_{n\to\infty} \ex(n,\cF)/n^r=0$
 --- and therefore their exact order of growth is not known.

\bsk



\bigskip

\nin \textbf{AMS 2000 Subject Classification:}
 05D05, 
 05C65, 
 68R05.

\end{abstract}

\bsk

\section{Introduction}

In this paper we study a Tur\'an-type problem on uniform hypergraphs,
 which is motivated by optimization of distributed data storage
 enabling secure data retrieval under a certain protocol.

\subsection{Terminology}

 \paragraph{Hypergraphs.}  A \emph{hypergraph} $H$ is a set system with vertex set
 $V(H)$ and edge set $E(H)$ where every edge $e \in E(H)$ is a
 nonempty subset of $V(H)$. The number of its vertices and edges is
 the \emph{order} and the \emph{size} of $H$, respectively.
  A hypergraph $H$ is called \emph{$r$-uniform}
  if  each edge of it contains precisely $r$ vertices.
  For short, sometimes we shall use the term   \emph{$r$-graph}
  for $r$-uniform hypergraphs.
  Graphs without loops
are just 2-uniform hypergraphs. A hypergraph $H_1$ is a
\emph{subhypergraph} of $H_2$ if $V(H_1) \subseteq V(H_2)$ and
$E(H_1) \subseteq E(H_2)$ holds, moreover we say that $H_1$ is an
\emph{induced subhypergraph} of $H_2$   if also $E(H_1)=\{e:\, e
\subseteq V(H_1) \enskip \wedge \enskip e \in E(H_2)\}$ holds. In
this paper graphs and hypergraphs are meant to be simple, that is
without loops and multiple edges, unless stated otherwise
explicitly.
  \bsk

  \paragraph{Tur\'an numbers.} Given hypergraphs $H$ and $F$, $H$ is said to be \emph{$F$-free} if $H$
 has no subhypergraph isomorphic to $F$. Similarly, if $\cF$ is a
 family of hypergraphs, $H$ is $\cF$-free if it  contains no
  subhypergraph isomorphic to any member of $\cF$.
  In the  problems considered here, the family $\cF$ contains
  $r$-graphs for a fixed $r \ge 2$ and the property to be
  $\cF$-free is considered only for $r$-graphs.

  In a Tur\'an-type (hypergraph) problem there is a given  collection $\cF$ of $r$-uniform
  hypergraphs   and
  the main goal is to determine or to estimate the
  \emph{Tur\'an number}
  $\ex (n, \cF)$ which is the maximum number of edges in an
  $\cF$-free   $r$-uniform hypergraph on $n$ vertices.
In 1941 Tur\'an \cite{Tu} determined
  $\ex (n, K_t)$, that is  the maximum size of a graph $G$ of order
  $n$ such that $G$ contains no complete subgraph on $t$ vertices.
  (The spacial case of $k=3$ was already solved in 1907 by Mantel
  \cite{Ma}.)
  Since then lots of famous results have been proved (see  the
  recent surveys \cite{FurS, Keev}), but many problems  especially
  among the ones concerning hypergraphs  seem  notoriously hard.
   \bsk

  \paragraph{Combinatorial batch codes.} The notion of batch code
    was introduced by Ishai, Kushilevitz, Ostrovsky
 and Sahai \cite{YKOS} to represent the distributed storage of
 $m$ items of data on   $n$ servers such that any at most $k$
  data items are recoverable  by submitting at most
 $t$   queries to each server.\footnote{In the main part of the literature
 notations $n$ and $m$ are used in reversed role. Here   the
 usual notation of hypergraph Tur\'an problems is applied for CBCs (as  done also in \cite{BB}).}
   In its combinatorial version \cite{PSW},
 `encoding' and `decoding' mean simply that the data items are stored on and read  from the
 servers. Its basic case, when the parameter $t$
 equals  1, can be defined as follows.\footnote{In this
 definition the vertices of the hypergraph represent the $n$ servers,
  the edges represent the $m$ data items, and an edge contains exactly
  those vertices which correspond to the servers storing the data items
  represented by the edge. Parameters $k$ and $t=1$ express the
  condition that every family of at most $k$ edges has a system of
  distinct representatives. Applying Hall's Theorem we obtain the
  definition in the form given here.}

  \tmz \item
  A \emph{combinatorial batch code (CBC-system)} with parameters $(m,k,n)$
  is a multihypergraph $H$ of order $n$ and size $m$, such that the union
   of any $i$ edges contains at least $i$ vertices for every $1 \le i \le k$.
   For given parameters $r,k,n$, satisfying   $r\ge 2$ and $r+1\le k \le n$,
   let $m(n,r,k)$ denote the maximum number $m$  of edges such
   that   an   $r$-uniform CBC-$(m,k,n)$-system exists.
  \etmz

  Optimization problems on combinatorial batch codes (mainly for the non-uniform case
  and under the condition $t=1$) were studied in
  \cite{BRR,BKMS,k=4,miamano,cbc-ENDM,cbc-AADM,PSW}.
  Recently,
  Balachandran and Bhattacharya \cite{BB} formulated the problem of
  determining the maximum size of  $r$-uniform CBC-systems as a Tur\'an multihypergraph
  problem. Clearly, an $r$-uniform multihypergraph $H$ is a CBC-system with parameter $k$ if
  and only if
  it has no subhypergraph  of order $i-1$ and size
  exactly $i$ for all $r+1 \le i \le k$.

  \bsk
  \paragraph{A problem of Brown, Erd\H{o}s and T.\ S\'os.}
  Brown, Erd\H{o}s and T.\ S\'os started to study the problems where, for
  fixed integers $2\le r \le v$ and $k \ge 2$,
  all $r$-graphs on $v$ vertices and with at least $k$ edges are
  forbidden to occur as a subhypergraph of an $r$-graph
   \cite{BES1}.\footnote{On graphs, the problem was first studied by
    Dirac in \cite{Di}.}
  The maximum
  size of such an $r$-graph  of order $n$ is denoted by
  $f^{(r)}(n,v,k)-1$. A general lower bound  on
  $f^{(r)}(n,v,k)$ was proved in \cite{BES1} and later further famous results
  were given for the cases $v\ge k$
  (see, e.g., \cite{RSz, EFR, SS1, SS2, AS}).
   In this paper, motivated by the optimization  problem on uniform CBCs,
    we will study a problem closely related to the   case $v \le k$.

   \bsk

   \paragraph{Our problem setting.}
   We shall consider Tur\'an-type problems for the following families of
   forbidden subhypergraphs.
    The upper index
 `$(r )$' in the notation indicates that the family
  consists of  $r$-graphs.
 \tmz
   \item $\cH^{(r)}(k,q)= \{H: \, |E(H)|-|V(H)| = q+1 \enskip \wedge \enskip    |E(H)| \le k
   \}$\\\\
    To study $\cH^{(r)}(k,q)$-free hypergraphs, we put the following restrictions on the parameters:
    \tmz
    \item[$\circ$] $r\ge 2$ \enskip (The problem would be trivial
    for the 1-uniform case.)
     \item[$\circ$] $k \ge q+r+1$ \enskip ($|E(H)| \le q+r$
     would imply
     $|V(H)| \le r-1$ and hence
        $\cH^{(r)}(k,q)=\emptyset$.)
       \item[$\circ$] $q \ge -r+1$ \enskip (Negative values can be
       allowed for $q$. But if $q\le -r$,   the family $\cH^{(r)}(k,q)$
          contains an $r$-graph with 1 edge and with at least $r$
       vertices, and hence $\ex (n, \cH^{(r)}(k,q))=0$ would follow.)
    \etmz
    \item $\cF^{(r)}(k,q)= \{H: \, |E(H)|-|V(H)| = q+1 \enskip \wedge \enskip    |E(H)|
    =    k\}$\\\\
     In general, $r\ge 2$, \enskip $k \ge q+r+1$ and $k \ge 2$ are
     assumed. Here we restrict ourselves to the  cases with $q \ge
     -r+1$.
     Note that $\cF^{(r)}(k,q)$ contains exactly those $r$-graphs which are
     forbidden in the Brown-Erd\H{o}s-S\'os problem with $v=k-q-1$,
          while $\cH^{(r)}(k,q)= \cup_{i=r+q+1}^k
     \cF^{(r)}(i,q)$.
   \etmz
   Moreover, for  $\cH^{(r)}(k,q)$ and
     $\cF^{(r)}(k,q)$,
       the family of multihypergraphs with the same defining
     property is denoted by $\cH^{(r)}_M(k,q)$ and
     $\cF^{(r)}_M(k,q)$, respectively.
     When the Tur\'an number relates  to the maximum size of a
     multihypergraph, the lower index $M$ is used, as well. For
     instance, $\ex_M(n,  \cH^{(r)}_M(k,q))$ denotes the maximum
     number of edges in a multihypergraph such that every $i$ edges
     cover at least $i-q-1$ vertices subject to  $q+r+1 \le i \le k$.
     Note that if $q=-r+1$,  already the presence of  edges with
     multiplicity 2
       is  forbidden   and consequently
     $\ex_M(n,\cH^{(r)}_M(k,-r+1))= \ex (n,\cH^{(r)} (k,-r+1))$.

     \bsk

     The next facts follow immediately from the definitions:
   \begin{eqnarray*}
   m(n,r,k) &=& \ex_M(n,\cH^{(r)}_M(k,0))\\
   \ex (n,\cH^{(r)} (k,q))  \le   \ex (n,\cF^{(r)} (k,q)) &=&  f^{(r)}(n,k-q-1,k)-1
   \le \ex_M (n,\cF^{(r)}_M (k,q))\\
   \ex(n,\cH^{(r)}(k,q)) &\le & \ex_M(n,\cH^{(r)}_M(k,q))
   \end{eqnarray*}

\subsection{Preliminaries and our results}

 The following general lower bound was proved
 by Brown, Erd\H{o}s and T.\ S\'os \cite{BES1} for
 $\cF^{(r)} (k,q)$-free $r$-graphs under the previously given
 conditions ($r \ge 2$, $k \ge q+r+1$ and $k \ge 2$).
   \begin{eqnarray}
       f^{(r)}(n,k-q-1,k)= \Omega( n^{r-1 + \frac{q+r}{k-1}}).
  \end{eqnarray}
 Paterson, Stinson and Wei \cite{PSW} proved that if $q=0$ but all the
 $r$-graphs from $\cH^{(r)}(k,0)$ are forbidden, the  lower
 bound (1) still remains valid\footnote{For the cases with $k-\lceil \log
 k \rceil \le r \le k-1$,   Balachandran and Bhattacharya \cite{BB}
   proved the better lower bound $m(n,r,k) = \Omega(n^r)$}:
   \begin{eqnarray*}
        m(n,r,k) \ge \ex(n, \cH^{(r)}(k,0))  =  \Omega(n^{r-1 + \frac{r}{k-1}}).
  \end{eqnarray*}
  We prove in Section 2 that  the lower bound (1) can be extended also to our general
  case:
  \begin{eqnarray}
       \ex(n, \cH^{(r)}(k,q))=\Omega( n^{r-1 + \frac{q+r}{k-1}}).
  \end{eqnarray}

   \bsk

  Concerning  upper bounds, our main result proved in Section 4 says that
\begin{eqnarray}
              \ex(n,\cH^{(r)}(k,q)) = \cO ( n^{r-1+ \frac{1}{\left \lfloor
 \frac{k}{q+r+1}\right\rfloor}})
       \end{eqnarray}
  for every fixed $r\ge 2$ and $k\ge  q+r+1$.
  The basis of the proof is $r=2$ (graphs), for which the order of
  the upper bound follows already from a theorem of Faudree and
  Simonovits \cite{FS}; in fact they only forbid a subfamily
  of $\cF^{(2)}(k,q)$.
   Under the stronger condition of excluding
  $\cH^{(2)}(k,q)$ instead of $\cF^{(2)}(k,q)$, however, a better and
  explicit constant can be derived on the former; and this can in turn
   be proved to be valid on the latter as well.
  For this reason, we do not simply derive the result from the one
   in \cite{FS} but prove the new upper bound in our Theorem \ref{graphs}.
   The more general result for
  hypergraphs is given in Theorem \ref{hgs}.
    In Section 4 we also prove that the same upper bound (3)
  is  valid for multihypergraphs, in fact not only the orders of these upper bounds are equal but
 also the relatively small leading coefficients are the same.

  Section 5 is devoted to exploring the connection between the Tur\'an numbers of
  $\cH^{(r)}(k,q)$ and $\cF^{(r)}(k,q)$.
    The general message there is that any later improvement in the estimates
   concerning $\cH^{(r)}(k,q)$ will automatically yield an improvement for
    $\cF^{(r)}(k,q)$ as well, and vice versa.\footnote{Obviously, by this
    principle, one should seek upper bounds for $\cH^{(r)}(k,q)$
     and lower bounds for $\cF^{(r)}(k,q)$.}
  By Theorem \ref{th-d}, if $r=2$ and the parameters  $k$ and $q$ are fixed, the difference
  is   bounded by a constant $d(k,q)$:
  $$f^{(2)}(n,k-q-1,k)-\ex(n,\cH^{(2)}(k,q)) \le d(k,q).$$
    For $r \ge 3$, by Theorem \ref{th-d-hg} we obtain the upper bound
   $$f^{(r)}(n,k-q-1,k)-\ex(n,\cH^{(r)}(k,q))=\cO(n^{r-1}),$$
   which is somewhat weaker but still strong enough to prove that
   the Tur\'an numbers $\ex(n,\cF^{(r)}(k,q))$
  and $\ex(n,\cH^{(r)}(k,q))$ have the same order of growth.
   On the other hand, the question of sharpness
  of Theorem \ref{th-d-hg} remains open:


  \bpm   \label{diff-inf}
   For the triplets\/ $(r,k,q)$ of integers in the range\/ $r\ge 2$,\/
    $q \ge -r+1$,\/ and\/ $k \ge q+r+1$, determine the infimum
     value\/ $s(r,k,q)$ of
    constants\/ $s\ge 0$ such that
      $$f^{(r)}(n,k-q-1,k)-\ex(n,\cH^{(r)}(k,q))=\cO(n^s)$$
  as\/ $n\to\infty$.
  \epm

  \bcj   \label{diff-min}
  The infimum\/ $s(r,k,q)$ in Problem \ref{diff-inf} is
   attained as minimum.
  \ecj

   Our Theorem \ref{th-d} shows that $s(2,k,q)=0$ holds
     for all pairs $(k,q)$ in the given range, and so
     Conjecture \ref{diff-min} is confirmed for $r=2$.

 \bsk

   At the end of this introductory section, we return to uniform combinatorial batch codes.
   The previous upper bound given for $m(n,r,k)$ in
   \cite{PSW} was improved recently by Balachandran and
   Bhattacharya
    \cite{BB}:
 \begin{eqnarray}
       m(n,r,k) &=&  \cO(  n^{r-  \frac{1}{2^{r-1}}}) \qquad\qquad  \mbox{if
       \enskip
       $3\le r\le k-1-\lceil \log k\rceil$}.
  \end{eqnarray}
  Our Corollary \ref{cbc}  yields a further  improvement  in the
  range $r \le k/2 -1$. Especially, we have
    \begin{eqnarray}
       m(n,r,k) = \cO\left( n^{r-1+ \frac{1}{\left\lfloor
 \frac{k}{r+1}\right\rfloor}}\right).
         \end{eqnarray}
      Comparing (4) and (5), the difference is significant already
      for  parameters complying with $3 \le r = k/2 -1$. For  these cases,
      (4) gives exponent $r-   1/2^{r-1}$ whilst
      our bound (5) yields exponent $r-1/2$.

 \section{Lower bound}

 In this section we prove a lower bound on $\exsk$
 whose order is the same as proved in \cite{BES1} for
 $f(n,k-q-1,k)$; that is, for the case when only the subhypergraphs
 on exactly $k-q-1$ vertices and with $k$ edges are forbidden.

\begin{Theorem} \label{lower}
 For all fixed triplets of integers\/ $r,k,q$
 with\/ $r \ge 2$,\/ $q \ge -r+1$ and\/ $k \ge r+q+1$ we have
 $$
   \exsk = \Omega(n^{r-1+\frac{q+r}{k-1}})
     = \Omega(n^{\frac{kr-k+q+1}{k-1}}).
 $$
\end{Theorem}

\pf
 We apply the probabilistic method. Our proof technique is similar
 to those  in \cite{BES1} and \cite{PSW}.
 We let $p=cn^{-1+\frac{q+r}{k-1}}$, where the constant
  $c=c(r,k,q)$ will be chosen later.
 Note that the lower bound $-r+1$ on $q$ implies
  $pn\ge cn^{\frac{1}{k-1}}$, i.e.\
  $pn$ tends to infinity with $n$ whenever $r,k,q$ are constants.

 Let $\hnp$ be the random $r$-uniform hypergraph of order $n$
  with edge probability $p$.
 That is, $\hnp$ has $n$ vertices, and for each $r$-tuple $S$ of
  vertices the probability that $S$ is an edge is $p$,
  independently of (any decisions on) the other $r$-tuples.
 We denote by $E$ the number of edges in $\hnp$, and by $F$
  the number of forbidden subhypergraphs in $\hnp$;
 by `forbidden' we mean that for some $i\le k$, some $i-q-1$
  vertices contain at least $i>0$ edges.

 We will estimate the expected value of $E-F$, more precisely
  our goal is to show that the inequality
   $\bbE(E-F)\ge \bbE(E)/2$ on the expected values is true for a
  suitable choice of the constant $c$.
 Once $\bbE(E-F)\ge \bbE(E)/2$ is ensured, we obtain that
  there exists a (non-random) hypergraph with twice as many
  edges as the number of its forbidden subhypergraphs,
  hence removing one edge from each of the latter we obtain
  a hypergraph with the required structure and with
  at least $\bbE(E)/2=\frac{p}{2} {n\choose r}$ edges.

 By the additivity of expectation we have
$$
  \bbE(E-F) = \bbE(E) - \bbE(F),
$$
  moreover it is clear by definition that
\begin{eqnarray} \label{E(E)}
  \bbE(E) = p\cdot\! {n\choose r}
    = ({\textstyle\frac{1}{r!}}+o(1))\cdot p\cdot n^r
      = ({\textstyle\frac{1}{r!}}+o(1))
        \cdot c \cdot n^{r-1+\frac{q+r}{k-1}}
\end{eqnarray}
  for any fixed $r$ as $n\to\infty$.
  Hence we need to find an upper bound on $\bbE(F)$.

 We consider the following set $I$ of those values of $i$ for which an
   $(i-q-1)$-element vertex
 subset is large enough to accommodate some forbidden subhypergraph:
 $$I=\left\{i : \enskip i \le {i-q-1\choose r} \quad \wedge \quad q+r+2 \le i \le
 k\right\}.$$
 It should be noted first that if $I=\emptyset$,
  then also $\cH^{(r)}(k,q)=\emptyset$ holds
 and  hence $\exsk= {n \choose r}$. In this case, the lower bound in the theorem is
 trivially valid, as the condition
 $k \ge r+q+1 \ge 2$ implies $(q+r)/(k-1) \le 1$.

 From now on, we assume that $I \neq \emptyset$.
 Consider any $i\in I$.
 On any $i-q-1$ vertices the number of ways we can select
  $i$ edges is ${{i-q-1\choose r}\choose i}$, and the
  probability for each of those selections to be a subhypergraph
  of $\hnp$ is exactly $p^{i}$.
 Since there are ${n\choose i-q-1}$ ways to select $i-q-1$
  vertices, we obtain the following upper bound:
\begin{eqnarray} \label{E(F)}
  \bbE(F) & \le &  \sum_{i\in I} {{i-q-1\choose r}\choose  i} \cdot
    p^{ i}\cdot\! {n\choose i-q-1} \nonumber \\
  & < & \sum_{i\in I} \frac{ {{i-q-1\choose r}\choose  i} }
    {(i-q-1)!}
    \cdot p^{ i}\, {n^{i-q-1}} \nonumber \\
  & < & {\textstyle \left( {\displaystyle \max_{i\in I} } \,
    \frac{ {{i-q-1\choose r}\choose i} } {(i-q-1)!}
      \right) } \cdot p^k \, n^{k-q-1} \cdot
    \sum_{i=q+r+2}^{k} (pn)^{i-k} \nonumber \\
  & \le & ( C_{k,q,r} + o(1) ) \cdot c^k \cdot
      n^{k-q-1-k\,(1- \frac{q+r}{k-1})}  \nonumber \\
  & = & ( C_{k,q,r} + o(1) ) \cdot c^k \cdot
      n^{r-1+\frac{q+r}{k-1}}
\end{eqnarray}
  where $C_{k,q,r}$ abbreviates the maximum value of
  $\frac{ {{i-q-1\choose r}\choose i} } {(i-q-1)!}$
  taken over the range $I$ of $i$.

 Compare the rightmost formula of (\ref{E(E)}) with
   (\ref{E(F)}).
 The terms in parentheses containing $o(1)$ are essentially
  constant, while the main part of (\ref{E(E)})
  is \mbox{$c \cdot n^{r-1+\frac{q+r}{k-1}}$} whereas that of
  (\ref{E(F)}) grows with $c^k \cdot n^{r-1+\frac{q+r}{k-1}}$.
Thus, choosing $c$ sufficiently small, the required
 inequality $\bbE(E-F)\ge \bbE(E)/2$ will hold for
 $n$ large.
 This completes the proof of the theorem. \qed
 \bsk

 \brm
 It can also be ensured (again by a suitable choice of $c$)
  that $\bbE(E-F)/\bbE(E)$ is arbitrarily close to 1.
 This is not needed for the proof above, but it may be of
  interest in the context of batch codes with specified rate
  (cf. e.g. \cite{YKOS}).
 \erm

\section{Upper bound for graphs}

  First we prove an upper bound on ex$(n,\cH^{(2)}(k,q))$.

 \thm \label{graphs}
 For every three integers\/ $q\ge -1$,\/ $k \ge 2q+6$ and\/ $n \ge
 k$, we have
    $$  \ex(n,\cH^{(2)}(k,q)) < C\cdot n^{1+ \frac{1}{\left \lfloor
 \frac{k}{q+3}\right\rfloor}}+(q+2)n,$$
 where\/ $C=(q+2)^{\frac{1}{\left \lfloor
 \frac{k}{q+3}\right\rfloor}}$.
 \ethm
 \pf Introduce the notation $h=\left \lfloor  \frac{k}{q+3}\right\rfloor$
 and assume for a contradiction that there exists a graph $G$ of order $n$
 in which, for every $q+3 \le i \le k$, every $i$ edges  cover at least $i-q$
 vertices and the number of edges in $G$ is
 $$ |E(G)|=m \ge C\cdot n^{1+ \frac{1}{h}}+(q +2)n.$$
 Thus, the average degree $\bar{d}(G)=\bar{d}$  satisfies
  $$ \bar{d}=\frac{2m}{n} \ge 2C \cdot  n^{\frac{1}{h}}+2(q+2).$$
 Moreover, every graph of average degree $\bar{d}$ has a subgraph of
 minimum degree greater than  $\bar{d}/2$.\footnote{Just
 delete sequentially the vertices of degree smaller than or equal to
 $\bar{d}/2$. After each single step the average degree is
 greater than or equal to $\bar{d}$. Hence, finally we obtain a
 subgraph of minimum degree greater than $\bar{d}/2$.}
  Hence, we have a subgraph $F$ with
 minimum degree $\delta(F)= \delta$ such that
 \begin{eqnarray}\delta > C \cdot n^{\frac{1}{h}}+q+2.
 \end{eqnarray}

 \msk
 \noindent\underline{{\it Claim A.}}\, The order of $F$ satisfies
 $$|V(F)| >     \frac{(\delta-q-2)^{h}}{q+2}  .$$
 \msk

 \noindent{\it Proof.} Choose a vertex $x$ of $F$ as a root and construct the
  breadth-first search tree (BFS-tree) of $F$ rooted in $x$.
  Let $L_i$ denote the set of vertices on the $i$th level of the
  BFS-tree, and introduce the notation
  $\ell_i=|L_i|$.
  The edges of $F$ not belonging to the BFS-tree will be called
  additional edges.

      First we consider the vertices of the first
      $h^*=\left \lfloor  \frac{k-q-1}{q+3}\right\rfloor$ levels
      and prove that each vertex $v \in L_i$
        is incident with at most $q+1$
      additional edges, if $0 \le  i \le h^*-1$.
      Assume to the contrary that there exist $q+2$ such additional
  edges and consider the union of paths on the BFS-tree connecting
  the end-vertices of these additional edges with the root vertex
  $x$.
  This means $q+3$ (not necessarily edge-disjoint) paths
  each of length at most $h^*$, and at least one of them (the path
  between $v$ and $x$) is of length at most $h^*-1$.
   They form a tree, let  the number of its edges
  be denoted by $p$.
  Together with the $q+2$ additional
  edges we have
  $$p+q+2 \le h^*-1+(q+2)h^*+q+2 = (q+3)h^* +q+1 \le k$$
   edges, which cover only $p+1$ vertices. This
  contradicts the assumed property of $G$. Therefore, we may have
  at most $q+1$ additional edges incident with  vertex
  $v$.

   Now, we prove a  bound on  the number $\ell_i$ of vertices on the $i$th
    level if $2 \le i \le h^*$. The sum of the vertex degrees
  over the set $L_{i-1}$ cannot be smaller than $ \delta \ell_{i-1}$.
    On the other hand, each of these $\ell_{i-1}$ vertices is incident with  at most
    $q+1$ additional edges, moreover there are
     $\ell_{i-1}+\ell_i$ edges of the BFS-tree each of them being  incident with exactly one vertex from
     $L_{i-1}$.
      As follows,
  \begin{eqnarray*}
  \delta \, \ell_{i-1} & \le & \ell_{i-1}+\ell_i+(q+1)\ell_{i-1}\\
  (\delta-q-2)\;  \ell_{i-1} &\le & \ell_i,
   \end{eqnarray*}
   for every $2 \le i \le h^*$.
   Since $\ell_1 \ge \delta -q-2$ is also true, the recursive
   formula gives
   \begin{eqnarray}
   |V(F)|\ge  \ell_{h^*} \ge (\delta-q-2)^{h^*} \ge \frac{(\delta-q-2)^{h^*}}{q+2}.
   \end{eqnarray}
   If $h=h^*$, that is if $k \equiv q+1$ or $q+2$ (mod \, $q+3$),
   this already proves Claim~A.

   In the other case we have $h=h^*+1$ and claim that every
   vertex $u\in L_{h-1}$ is incident with at most $q+1$
   additional edges whose other end is in $L_{h-2} \cup L_{h-1}$.
    Then, assume for a contradiction that there are at least $q+2$ such edges.
   Again, take these $q+2$ additional edges together with
   the paths in the BFS-tree connecting
  their ends with the root.
  In this subgraph we have only at most
  $(q+3)(h-1) +q+2 <k$ edges,
   which
   cover fewer vertices by $q+1$ than the number of edges. Proved by this
   contradiction, we have at most $q+1$ additional edges of the
   described type.

   A similar argumentation  shows that each $w \in L_h$ might be incident with at most $q+1$
    additional edges whose other end is in $L_{h-1}$.
  Assuming the presence of $q+2$ such edges, we have at most
   $h+ (q+2)(h-1)+q+2 \le k$ edges  together with the
   paths between their ends and the root. Moreover,
   this cardinality exceeds the number of covered vertices by $q+1$.  Thus, we have a
   contradiction, which proves the property stated  for $w$.

   By these two bounds on the number of additional edges we can
   estimate the sum $s$ of vertex degrees over $L_{h-1}$ as follows:
   $$ \delta\, \ell_{h-1} \le s \le \ell_{h-1}+\ell_h +
   (q+1)\ell_{h-1} + (q+1) \ell_h.$$
   Together with (9) this implies
    $$ |V(F)|\ge \ell_h \ge \frac{\delta-q-2}{q+2}\; \ell_{h-1} \ge
      \frac{(\delta-q-2)^{h}}{q+2},$$
   and proves    Claim A. \dia

  \bsk

  Turning to  graph $G$,
  inequality (8) and Claim A yield the contradiction
  $$
  n  \ge |V(F)| > \frac{\left(C \cdot n^{1/h}\right)^h}{q+2}=n.
 $$
  Therefore, in a $\cH^{(2)}(k,q)$-free graph  the number of edges must be smaller than
  $C\cdot n^{1+ 1/h}+(q+2)n$,
  as stated in the theorem. \qed

  \bc
   \label{multi-graphs}
  For every three integers\/ $q\ge -1$,\/ $k \ge 2q+6$ and\/ $n \ge
  k$,  we have
    $$  \ex_M(n,\cH^{(2)}_M(k,q)) < C\cdot n^{1+ \frac{1}{\left \lfloor
 \frac{k}{q+3}\right\rfloor}}+(q+2)n,$$
 where\/ $C=(q+2)^{\frac{1}{\left \lfloor
 \frac{k}{q+3}\right\rfloor}}$.
  \ec
  \pf The BFS-tree of a multigraph $G$ is meant as a simple graph. That
  is, if an edge $uv$ has multiplicity $\mu >1$ in $G$, and $uv$ is an
  edge in the BFS-tree, then only one edge $uv$ belongs to the tree,
  the remaining $\mu -1$ copies are   additional edges.
  With this setting every detail of the previous proof remains valid for
  multigraphs. \qed

\section{Upper bound for hypergraphs}

  In this section we study the  problem for hypergraphs. The upper bound on
 $\ex(n,\cH^{(r)}(k,q))$ will be obtained by using  Theorem
\ref{graphs}.

  \thm \label{hgs}
  Let\/ $n$, $k$, $r$  and\/ $q$ be integers such that\/ $r\ge 2$, $q\ge -r+1$ and $n\ge k \ge
  2q+2r+2$, moreover let $C'=(q+r)^{\frac{1}{\left \lfloor
 \frac{k}{q+r+1}\right\rfloor}}$.
  Then,
  $$   \exsk  < \frac{2C'}{r!}\cdot n^{r-1+ \frac{1}{\left \lfloor \frac{k}{q+r+1}\right
   \rfloor}}+\frac{2(q+r)}{r!}\cdot n^{r-1}. $$
  \ethm

  \pf Consider an $\cH^{(r)}(k,q)$-free $r$-graph $H$. Let  its
  order and size be denoted by $n$ and $m$, respectively.
  For a set $S\subseteq V(H)$ denote by $d(S)$ the number of edges
  of $H$ which contain $S$ entirely. By double counting we have
  $$\sum_{S \subset V(H), \enskip |S|=r-2}d(S)= m\;{r \choose r-2},$$
    and for the average value $\bar{d}_{r-2}$ of $d(S)$ over the
  $(r-2)$-element subsets of $V(H)$
  $$
  \bar{d}_{r-2}= m\;\frac{{r \choose r-2}}{{n \choose r-2}}
  $$
  holds. Thus, there exists an $S^* \subset V(H)$ of cardinality $r-2$
  satisfying
  $$
   d(S^*) \ge m\;\frac{{r \choose r-2}}{{n \choose r-2}}.
  $$

    Deleting the edges which do not contain $S^*$ entirely,
   in addition deleting  the $r-2$ vertices of $S^*$ from the remaining edges,
   we obtain a graph $G$ with $V(G)=V(H)$ and
  $$E(G)= \{e \setminus S^*: S^* \subset e \enskip \wedge \enskip e \in E(H)\}, \quad
    \quad |E(G)| \ge m\;\frac{{r \choose r-2}}{{n \choose r-2}} .$$
   Since  every $i$ edges ($i \le k$) cover at least $i-q$ vertices in $H$,
      every $i$ edges cover at least $i-q-r+2$ vertices
   in $G$.
   Moreover, the conditions given in Theorem \ref{graphs} hold for
   $n'=n$, $k'=k$ and $q'=q+r-2$. Then, we obtain
   \begin{eqnarray}
   m\;\frac{{r \choose r-2}}{{n \choose r-2}}
   \le |E(G)| <
   (q+r)^{\frac{1}{\left \lfloor \frac{k}{q+r+1}\right   \rfloor}}
    n^{1+ \frac{1}{\left \lfloor \frac{k}{q+r+1}\right   \rfloor}}
    +(q+r)n,
    \end{eqnarray}
    from which
    $$
    m < \frac{2C'}{r!}\; n^{r-1+ \frac{1}{\left \lfloor \frac{k}{q+r+1}\right
   \rfloor}}+\frac{2(q+r)}{r!}\cdot n^{r-1}
    $$
    follows. This implies the same upper bound for $\ex(n,\cH^{(r)}(k,q))$.
    \qed

   \bsk

   The above proof remains valid if  the $r$-graph $H$  is allowed
   to have multiple edges. The only difference is that
   we must refer to Corollary \ref{multi-graphs} instead of Theorem
   \ref{graphs}.
    Hence, for
   multihypergraphs the same upper bound can be stated. In addition,
   since $m(n,r,k)=\ex_M(n,\cH^{(r)}_M(k,0))$, we obtain a new upper bound
   for the maximum size $m(n,r,k)$ of $r$-uniform CBC-systems with parameters
   $n$ and $k$.
   \bc \label{multi-hgs}
     Let\/ $n$, $k$, $r$  and\/ $q$ be integers such that\/ $r\ge 2$, $q\ge -r+1$ and $n\ge k \ge
  2q+2r+2$, moreover let $C'=(q+r)^{\frac{1}{\left \lfloor
 \frac{k}{q+r+1}\right\rfloor}}$.
  Then,
  $$   \ex_M(n,\cH^{(r)}_M(k,q))  < \frac{2C'}{r!}\cdot n^{r-1+ \frac{1}{\left \lfloor \frac{k}{q+r+1}\right
   \rfloor}}+\frac{2(q+r)}{r!}\cdot n^{r-1}. $$
   \ec

   \bc \label{cbc}
     Let\/ $n$, $k$, $r$   be integers such that\/ $r\ge 2$   and $n\ge k \ge
  2r+2$, moreover let $C''= r^{\frac{1}{\left \lfloor
 \frac{k}{r+1}\right\rfloor}}$.
  Then,
  $$   m(n,r,k)  < \frac{2C''}{r!}\cdot n^{r-1+ \frac{1}{\left \lfloor \frac{k}{r+1}\right
   \rfloor}}+\frac{2}{(r-1)!}\cdot n^{r-1}. $$
   \ec

 \section{Asymptotic equality of Tur\'an numbers}

  Up to this point we were concerned with the problem of $\cH^{(r)}(k,q)$-free
 hypergraphs; it is different from the one studied by
 Brown, Erd\H{o}s and T.\ S\'os  \cite{BES2,BES1}, where only the
 subhypergraphs with exactly $k-q-1$ vertices and $k$ edges are
 forbidden.
  In this section we show that
    $\ex (n,\cH^{(r)}(k,q))$  and $f^{(r)}(n,k-q-1,k)-1$
 are asymptotically equal.
  For graphs ($r=2$), our result is better as   there exists a constant upper bound
 (depending only on $k$ and $q$) on their difference.
   As a consequence, we   obtain a new upper bound on
   $f^{(2)}(n,v,k)$ subject to   $ v \ge (k+4)/2$.

  First we prove the following lemma. For fixed parameters $k$, $q$ and
  for a given graph $G$, a subgraph $G'$ is said to be  \emph{forbidden (for $(k,q)$)} if
  $G' \in  \cH^{(2)}(k,q)$, moreover $G'$ is \emph{maximal forbidden (for $(k,q)$)}, if it
  cannot be extended into a  forbidden subgraph of larger order.
\bl \label{lemma-d}
   Let $k$ and $q$ be integers such that $q \ge -1$ and $k \ge q+3$,
   and let $G$ be a graph of order at least $k-q-1$. If a subgraph $G'\subset G$ is
   maximal forbidden for $(k,q)$, then either $G'$ has $k$ edges or
   it is the union of one or more components of $G$.
\el
 \pf
 Assume that $G'$ is a forbidden subgraph of $G$ and
 $|E(G')| < k$.
   If there exists an edge $uv \in E(G)$ such
   that $u \in V(G')$ and $v \in V(G) \setminus V(G')$, then the
   subgraph $G''$ obtained by extending $G'$ with the vertex $v$ and
   with the edge $uv$ satisfies $|E(G'')|-|V(G'')|=q+1$ and
   $|E(G'')|= |E(G')|+1 \le k$. Hence $G''$ is   forbidden for
   $(k,q)$
   and consequently, $G'$ is not
   maximal forbidden.
   On the other hand, if the subgraph of $G$ which is induced by $V(G')$  contains some
   edge $e$ not in $G'$, then with any vertex $v \in V(G) \setminus
   V(G')$, the subgraph $G' + e +v$ is forbidden for $(k,q)$ and
   again, $G'$ is not a maximal forbidden subgraph.
   Therefore, if $G'$  is of order smaller than $k$ and
   it is a maximal forbidden subgraph for $(k,q)$, then $G'$
   is a component of $G$, or it is the union of some components of
   $G$.
   \qed
   \bsk

   Clearly,  $f^{(2)}(n,k-q-1,k) \ge \ex (n, \cH^{(2)}(k,q))$.
       The
   following theorem states that the difference between them is
   bounded by a constant, once the parameters $k$ and $q$ are fixed.

   \thm \label{th-d}
   For every pair $k,q$ of integers satisfying  $q \ge -1$ and $k \ge
   q+3$ there exists a constant $d=d(k,q)$ such that for every $n
   \ge k-q-1$,
   $$f^{(2)}(n,k-q-1,k) - \ex (n, \cH^{(2)}(k,q)) \le d.$$
   \ethm
   \pf For given parameters $k$ and $q$ first define
   $z:= \min \{i: \enskip q+3 \le i \le {i-q-1 \choose 2} \}$.
   If $k >z$, there is no forbidden subgraph for $(k,q)$ and
   consequently,
    $f^{(2)}(n,k-q-1,k) = \ex (n, \cH^{(2)}(k,q))= {n \choose 2}$.
    Otherwise, $z$ is
    the possible minimum size of a subgraph forbidden for $(k,q)$.
   By Theorem \ref{lower}
   $$\ex (n,   \cH^{(2)}(k,q))=\Omega(n^{1+\frac{q+2}{k-1}})$$
    holds, thus
   there exists an $n_0$ (depending only on $k$ and $q$) such that
   for all $n\ge n_0$
   \begin{eqnarray*} \label{1-5fej}
   \frac{z}{z-q-1}\cdot n \le \ex (n,   \cH^{(2)}(k,q)).
   \end{eqnarray*}
   Consequently, the following finite  maximum exists:  \label{2-5fej}
    \begin{eqnarray}
     d= \max \left(\left\{\frac{z}{z-q-1}\cdot n - \ex (n,   \cH^{(2)}(k,q))+1:
     \enskip n \in {\mathbb{N}}\right\}\cup \{1\}\right) .
    \end{eqnarray}
   We   claim that $d$ is a suitable constant
   for our theorem. To prove this, let us consider an $\cF(k,q)$-free graph $G$
    on $n$ vertices and with $f^{(2)}(n,k-q-1,k)-1$ edges. If $G$ is
   $\cH^{(2)}(k,q)$-free as well, $f^{(2)}(n,k-q-1,k)-1$ is equal to $\ex (n,
   \cH^{(2)}(k,q))$, and since $d \ge 1$, the theorem holds for $k$, $q$
   and $n$.

     In the other case, $G$ contains a  subgraph  $G_1$ maximal
     forbidden for $(k,q)$. Clearly, $G_1$ has fewer than $k$ edges,
     hence by Lemma \ref{lemma-d}, $G_1$ is an induced subgraph and there is no edge
     between $V(G_1)$ and $V(G)\setminus V(G_1)$.
     Then, the remaining subgraph $G-G_1$ is either $\cH^{(2)}(k,q)$-free or
     contains a subgraph  $G_2$ of size smaller than $k$, which is maximal forbidden for
     $(k,q)$.
     Iteratively applying this procedure, finally we have vertex-disjoint maximal
     forbidden subgraphs $G_1, \dots G_j$ and the $\cH^{(2)}(k,q)$-free
     subgraph $G'$ induced by $V(G)\setminus \cup_{i=1}^j V(G_i)$, such that each edge of $G$
     is contained in exactly one of $G',G_1, \dots G_j$.
      As $q+1 \ge 0$ and for every $1\le i \le j$  we have $z \le |E(G_i)|\le
      k-1$,  applying  Lemma \ref{lemma-d}, we obtain
      $$
     \frac{|E(G_i)|}{|V(G_i)|}= \frac{|E(G_i)|}{|E(G_i)-q-1|} \le
     \frac{z}{z-q-1}.
     $$
     Using   notations $n_1=\sum_{i=1}^j |V(G_i)|$ and
     $n_2=|V(G')|=n-n_1$, moreover the definition (11) of $d$
    \begin{eqnarray*}
    |E(G)|&=& f^{(2)}(n,k-q-1,k)-1  \le \frac{z}{z-q-1}\cdot n_1 +
     \ex (n_2,   \cH^{(2)}(k,q))\\
     &\le& \ex (n_1,   \cH^{(2)}(k,q))+d-1 + \ex (n_2,   \cH^{(2)}(k,q))\\
     &\le& \ex (n,   \cH^{(2)}(k,q))+d-1,
     \end{eqnarray*}
     which yields
     $$
   f^{(2)}(n,k-q-1,k)- \ex (n,   \cH^{(2)}(k,q)) \le d,
   $$
   as stated. \qed
   \bsk


 \bc
  Let\/ $v$ and $k$   be integers such that\/
 $2 \le v \le  k  $  and let\/ $C=(k-v+1)^{\frac{1}{\left \lfloor
 \frac{k}{k-v+2}\right\rfloor}}$.  Then, there exists a constant $D$
 such that for every $n$
   $$  f^{(2)}(n,v,k) < C\cdot n^{1+ \frac{1}{\left \lfloor
 \frac{k}{k-v+2}\right\rfloor}}+(k-v+1)n +D.$$
  \ec

\pf  Let $q$ denote $k-v-1$.  Then, under  the given conditions
  we have
 $-1 \le q \le k-3 $  and $C=(q+2)^{\frac{1}{\left \lfloor
 \frac{k}{q+3}\right\rfloor}}$.
 Theorems \ref{graphs} and \ref{th-d} immediately imply the
 existence of a constant $D$ such that for every $n$
    $$  f^{(2)}(n,k-q-1,k) < C\cdot n^{1+ \frac{1}{\left \lfloor
 \frac{k}{q+3}\right\rfloor}}+(q+2)n +D.$$
 This is equivalent to the statement of the corollary. \qed
 \bsk

 \thm \label{th-d-hg}
   For every four integers\/ $r,k,q$ and $n$   satisfying\/  $r \ge 2$ and\/ $2 \le q +r+1 \le
   k \le n$,
      $$f^{(r)}(n,k-q-1,k) - \ex (n, \cH^{(r)}(k,q)) \le (k-1) {n-1\choose r-1}$$
      holds. Hence, for every fixed\/ $r$, $k$, and\/ $q$ we have
      $$f^{(r)}(n,k-q-1,k) = (1+o(1)) \, \ex (n, \cH^{(r)}(k,q)).$$
   \ethm

 \bsk

\pf
Consider any extremal $r$-graph $H^*$ for $\cF^{(r)}(k,q)$
 on the $n$-element vertex set $V$.
By definition, $H^*$ is $\cF^{(r)}(k,q)$-free. If $H^*$ is also
$\cH^{(r)}(k,q)$-free, then
 $f^{(r)}(n,k-q-1,k) = \ex (n, \cH^{(r)}(k,q))$ holds
 and we have nothing to prove.
Otherwise we select the longest possible sequence of
 subhypergraphs $H_i\subset H^*$ ($i=1,2,\dots,\ell$)
  under the following conditions:
 \begin{itemize}
  \item
   Each $H_i$ is isomorphic to some member of
   $\cH^{(r)}(k,q) \smin \cF^{(r)}(k,q)$.
  \item
   Under the previous condition,
   $H_1$ is maximal in $H^*$.
  \item
   Under the previous conditions,
   $H_i$ is maximal in $H^*\smin\bigcup_{j=1}^{i-1} H_j$
   for each $2\le i\le\ell$.
 \end{itemize}
Eventually we obtain an $\cH^{(r)}(k,q)$-free hypergraph from $H^*$
 by removing at most $(k-1)\cdot\ell$ edges, because
 each $H_i$ has at most $k-1$ edges.
Thus, the proof will be done if we prove that
 $\ell\le{n-1\choose r-1}$ holds.

Let $e_i$ be an arbitrarily chosen edge of $H_i$ and let $f_i$ be
 an $(r-1)$-element subset of $e_i$, which we fix (again arbitrarily)
 for $i=1,2,\dots,\ell$.
Should $f_i\subset e_j$ hold for some $1\le i< j\le\ell$, the
hypergraph
 $H_i\cup\{e_j\}$ would also be isomorphic to some member of $\cH^{(r)}(k,q)$.
This contradicts the choice (maximality) of $H_i$.
 Consequently, for all $i=1,2,\dots,\ell$ we have:
 \begin{itemize}
  \item $|f_i|=r-1$,
  \item $|V\smin e_i|=n-r$,
  \item $f_i\cap (V\smin e_i)=\emptyset$,
  \item $f_i\cap (V\smin e_j)\ne\emptyset$ whenever $1\le i< j\le\ell$.
 \end{itemize}
Thus, applying a theorem of Frankl \cite{Fr},\footnote{Set pairs
  with prescribed intersection properties can be applied
 in many kinds of extremal problems (not only on graphs and hypergraphs).
A detailed account on those methods and results is given in the
 two-part survey \cite{T1,T2}.}
 the number of set pairs
 $(f_i, V\smin e_i)$ is at most ${(r-1)+(n-r)\choose r-1}
  = {n-1\choose r-1}$.
\qed

\bsk

  \bc
  Let\/ $r$, $v$, $k$   be integers such that\/ $r\ge 2$ and
 $(k+2r)/2 \le v \le  k+r-2  $  and let\/ $C=(k+r-v-1)^{\frac{1}{\left \lfloor
 \frac{k}{k+r-v}\right\rfloor}}$. Then,
 $$ f^{(r)}(n,v,k) \le \frac{2C}{r!}\cdot n^{r-1+\frac{1}{\left \lfloor
 \frac{k}{k+r-v}\right\rfloor}}+ \cO(n^{r-1}).$$
 \ec

{}


\begin{thebibliography}{}

\bibitem{AS}
 N. Alon and A. Shapira, \emph{On an Extremal Hypergraph Problem of Brown, Erd\H{o}s and   S\'os},
 Combinatorica 26 (2006), 627--645.


 \bibitem{BB}
  N. Balachandran and S. Bhattacharya,
  \emph{On an Extremal Hypergraph Problem Related to Combinatorial
  Batch Codes}, Manuscript (2012).

 \bibitem{BRR}
  S. Bhattacharya, S. Ruj and B. Roy,
 \emph{Combinatorial batch codes: A Lower Bound and Optimal
Constructions}, Advances in Mathematics of Communications 6  (2012),
165--174.

\bibitem{BES2}
W. G. Brown, P. Erd\H{o}s and V. T. S\'os, \emph{On the existence of
triangulated spheres in 3-graphs and related problems}, Periodica
Mathematica Hungarica, 3 (1973),   221--228.


\bibitem{BES1}
W. G. Brown, P. Erd\H{o}s and V. T. S\'os, \emph{Some extremal
problems on r-graphs}, New directions in the theory of graphs (Proc.
Third Ann Arbor Conf., Univ. Michigan, Ann Arbor, Mich., 1971),
53--63, Academic Press, New York, 1973.

\bibitem{BKMS}
   R. A. Brualdi, K. P. Kiernan, S. A. Meyer and M. W.
 Schroeder,
  \emph{Combinatorial batch codes and transversal matroids},
 Adv. Math. Commun., 4 (2010), 419--431.
  Erratum \textit{ibid.} p.~597.

  \bibitem{k=4}
  Cs. Bujt\'as and Zs. Tuza,
  \emph{Optimal batch codes: Many items or
   low retrieval requirement},
  Adv. Math. Commun.,  {\bf 5}  (2011), 529--541.

 \bibitem{miamano}
  Cs. Bujt\'as and Zs. Tuza,
  \emph{Optimal combinatorial batch codes
 derived from dual systems},
 Miskolc Math. Notes 12 (1) (2011), 11--23.



\bibitem{cbc-ENDM}
  Cs. Bujt\'as and Zs. Tuza,
 \emph{Combinatorial batch codes:
 Extremal problems under Hall-type conditions},
 Electr. Notes in Discrete Math. 38 (2011), 201--206.

\bibitem{cbc-AADM}
   Cs.\ Bujt\'as and Zs.\ Tuza:
   \emph{Relaxations of Hall's Condition: Optimal batch codes with multiple
   queries},
    Applicable Analysis and Discrete Mathematics, 6 (1) (2012),
   72--81.

\bibitem{Di}
  G. Dirac:
  \emph{Extensions of Tur\'an's theorem on graphs},
   Acta Math. Acad. Sci. Hungar. 14 (1963), 417--422.

\bibitem{EFR}
 P. Erd\H{o}s, P. Frankl and V. R\"odl,
 \emph{The asymptotic number of graphs not containing a
fixed subgraph and a problem for hypergraphs having no exponent},
Graphs Combin. 2 (1986), 113--121.

\bibitem{FS}
R. Faudree and M. Simonovits,
 \emph{On a class of degenerate extremal graph problems},
 Combinatorica 3 (1983), 83--93.

\bibitem{Fr}
P. Frankl, \emph{An extremal problem for two families of sets},
 European J. Combin. 3 (1982), 125--127.

\bibitem{FurS}
Z. F\"uredi and M. Simonovits,
 \emph{The history of degenerate (bipartite) extremal graph problems},
 in: Erd\H os Centennial (L. Lov\'asz et al., Eds.),
Bolyai Society Mathematical Studies 25 (2013), 169--264.

\bibitem{GST}
J. R. Griggs, M. Simonovits and G. R. Thomas,
 \emph{Extremal graphs with bounded densities of small subgraphs},
 Journal of Graph Theory,
   29  (1998), 185--207.

\bibitem{YKOS}
  Y. Ishai, E. Kushilevitz, R. Ostrovsky and A. Sahai,
  \emph{Batch codes and their applications},
 In: Proceedings of the 36th Annual ACM Symposium on Theory of
 Computing, ACM Press, New York, 2004, 262--271.

\bibitem{Keev}
P. Keevash, \emph{Hypergraph Turan problems}, Surveys in
Combinatorics, Cambridge University Press, 2011, 83--140.

\bibitem{Ma}
 W. Mantel, \emph{Problem 28}, Wiskundige Opgaven 10 (1907), 60--61.

\bibitem{PSW}
  M. B. Paterson, D. R. Stinson and R. Wei,
 \emph{Combinatorial batch codes},
 Adv. Math. Commun., {\bf 3} (2009), 13--27.

\bibitem{RSz}
 I. Z. Ruzsa and E. Szemer\'edi,
\emph{Triple systems with no six points carrying three triangles},
Combinatorics (Proc. Fifth Hungarian Colloq., Keszthely, 1976), Vol.
II, pp. 939--945, Colloq. Math. Soc. J ´nos Bolyai 18,
North-Holland, Amsterdam-New York, 1978.

\bibitem{SS1}
 G. N. S\'ark\"ozy and S. M. Selkow, \emph{An extension of the
 Ruzsa-Szemer\'edi Theorem}, Combinatorica 25 (2005), 77--84.

\bibitem{SS2}
 G. N. S\'ark\"ozy and S. M. Selkow, \emph{On a Tur\'an-type hypergraph problem
of Brown, Erd\H{o}s and T. S\'os}, Discrete Math. 297 (2005), 190--195.

\bibitem{Tu}
P. Tur\'an, \emph{On an extremal problem in graph theory} (in
Hungarian), Mat. Fiz. Lapok 48 (1941), 436--452.

\bibitem{T1}
 Zs. Tuza,
Applications of the set-pair method in extremal hypergraph theory.
    ``\,Extremal Problems for Finite Sets\,''
    (P. Frankl et al., eds.),
Bolyai Society Mathematical Studies 3,
     1994, 479--514.

\bibitem{T2}
 Zs. Tuza,
Applications of the set-pair method in extremal problems, II.,
``\,Combinatorics, Paul Erd\H os is Eighty\,''
(D. Mikl\'os et al., eds.),
Bolyai Society Mathematical Studies 2,
   1996, 459--490.


\end{thebibliography}
\end{document}